\documentclass{psp}
\usepackage{amssymb,indentfirst,amsfonts,amsmath,graphicx}
\usepackage{amscd}
\usepackage{graphicx,psfrag}

\DeclareMathAlphabet{\mathsfsl}{OT1}{cmss}{m}{sl}
\newcommand{\tensor}[1]{\mathsfsl{#1}}
\newcommand{\homo}{{\mathrm H}}

\newtheorem{theorem}{Theorem}[section]
\newtheorem{corollary}[theorem]{Corollary}
\newtheorem{proposition}[theorem]{Proposition}
\newtheorem{lemma}[theorem]{Lemma}

\newnumbered{definition}[theorem]{Definition}
\newnumbered{notation}[theorem]{Notation}

\newnumbered{remark}[theorem]{Remark}

\begin{document}

\title[Pseudo-Anosov extensions and degree one maps]
  {Pseudo-Anosov extensions and \\degree one maps between hyperbolic surface bundles}

\author[Michel Boileau, Yi Ni and Shicheng Wang]{MICHEL BOILEAU\\
Laboratoire \'Emile Picard,\addressbreak Universit\'e Paul
Sabatier, TOULOUSE Cedex 4, France\addressbreak e-mail\textup{:
\texttt{boileau@picard.ups-tlse.fr}}\nextauthor
YI NI\\
Department of Mathematics,\addressbreak Princeton University,  NJ
08544, U. S. A.\addressbreak e-mail\textup{:
\texttt{yni@math.princeton.edu}}\and\
SHICHENG WANG\\
LMAM, Department of Mathematics, \addressbreak Peking University,
Beijing 100871 China\addressbreak e-mail\textup{:
\texttt{wangsc@math.pku.edu.cn}} }

\maketitle

\begin{abstract}
 Let $F',F$ be any two closed orientable surfaces of genus $g'>g\ge 1$,
and $f:F\to F$ be any pseudo-Anosov map. Then we can ``extend" $f$
to be a pseudo-Anosov map $f':F'\to F'$ so that there is a fiber
preserving degree one map $M(F',f')\to M(F,f)$ between the
hyperbolic surface bundles. Moreover the extension $f'$ can be
chosen so that the surface bundles $M(F',f')$ and $M(F,f)$ have
the same first Betti numbers.
\end{abstract}

\section{Introduction}

All surfaces are oriented and all automorphisms on surfaces are
orientation preserving.

Let $F$ be an oriented closed surface of genus $g\ge 1$, and
$f:F\to F$ be an automorphism. We denote the surface bundle with
fiber $F$ and monodromy $f$ by $M(F,f)$.

\begin{definition}\label{defPA}
Suppose $G$ is a compact surface of genus $g\ge 1$. A circle $c$
on $G$ is {\it essential} if $c$ is neither contractible nor
boundary parallel. An automorphism $f$ of $G$ is {\it
pseudo-Anosov} if $f^n(c)$ is not isotopic to $c$ for any
essential circle $c\subset G$ and any integer $n$. (Note in the
case $G$ is a torus, the term ``pseudo-Anosov" we define here is
usually known as ``Anosov".)
\end{definition}

\begin{remark}Our definition of pseudo-Anosov
maps is slightly different from the more standard definition in
the literature. Pseudo-Anosov maps in our sense should be
considered as ``maps isotopic to a pseudo-Anosov map" in the
standard sense.
\end{remark}

Profound theories of Nielsen-Thurston and of Thurston in 2- and
3-dimensional topology tell us that pseudo-Anosov is the most
important class of surface automorphisms, and when $\chi(F)<0$,
$M(F, f)$ is a hyperbolic 3-manifold if and only if $f$ is
pseudo-Anosov.

\begin{theorem}\label{mainthm}
Let $F_s,F_t$ be closed orientable surfaces of genus $g_s$, $g_t$
respectively, $g_s>g_t\ge1$, and  $f_t:F_t\to F_t$ be a
pseudo-Anosov map. Then

(1) There exists a hyperbolic 3-manifold $M(F_s,f_s)$, such that
the is a fiber preserving degree one map $P: M(F_s,f_s)\to
M(F_t,f_t)$. (Here the subscript $s$ means ``source", and the
subscript $t$ means ``target".)

(2) Moreover the $f_s$ in (1) can be chosen so that $M(F_s,f_s)$
and $M(F_t,f_t)$ have the same first Betti numbers.
\end{theorem}

Motivation for Theorem \ref{mainthm} is from \cite{BW}, where the
following facts were proved:

(1) For each 3-manifold $M$, there is a degree one map
$f:M(F_s,f_s)\to M$ such that $M(F_s,f_s)$ is hyperbolic and
$\beta_1(M(F_s,f_s))=\beta_1(M)+1$.

(2) If there is a degree one map $f: M(F_s, f_s)\to M$ with
$\beta_1(M(F_s, f_s))=\beta_1(M)$ and $M$ is irreducible, then $M$
is a surface bundle and $f$ can be homotoped to a fiber preserving
one.

It is natural to wonder how to find fiber preserving degree one
maps between non-homeomorphic hyperbolic surface bundles (of the
same first Betti numbers). In Section~2, we will prove Theorem
\ref{mainthm} (1). In Section 3, by modifying the proof in
Section~2, we will prove Theorem \ref{mainthm} (2).

The proof of Theorem \ref{mainthm} relies on an extension process
from the pseudo-Anosov map $f_t$ on $F_t$ to a pseudo-Anosov map
$f_s$ on $F_s$, which is delicate and somewhat complicated.

We will outline  this process, i.e., for given $M(F_t, f_t)$ and
$F_s$, how to find $f_s$. In this outline we assume that $g_t\ge
2$. This process in Section 2 is divided into three steps.

Step 1. Fix a disk  $D\subset F_t$ and let
$V=F_t-\textrm{int}(D)$. We can assume that $f_t|D=\textrm{id}|D$
up to isotopy. Then as a restriction of a pseudo-Anosov map,
$f_t|V:V\to V$ is a pseudo-Anosov map (Lemma \ref{ResPA}).

Step 2. We will construct two embedding $e_0, e_1: V\to F_s$ such
that (1) $e_0(\partial V)$ and $e_1(\partial V)$ are not homotopic
in $F_s$,  (2) two pinches $p_0, p_1 : F_s\to F_t$ (see Definition
\ref{pinchdef}) defined by  $p_j\circ e_j = \text{id}_V :V\to V$
are homotopic (Lemma \ref{3exsitpinch}).

Step 3. The two embeddings $e_1$ and $e_2$ in step 2 also provided
a homeomorphism $\bar f_t := e_1 \circ f_t| \circ e_0^{-1}:
e_0(V)\to e_1(V)$. With properties of  $e_1$ and $e_2$ described
in Step 2, we will be able to extend $\bar f_t$ to a pseudo-Anosov
map $f_s: F_s\to F_s$ (Proposition \ref{4PAextension}).

Then clearly $p_1\circ f_s = f_t\circ p_0$, hence there exists a
fiber preserving degree one map $P: M(F_s, f_s)\to M(F_t, f_t)$
(Lemma \ref{1fiberprsv}). This finishes the proof of Theorem
\ref{mainthm} (1).

Now we give more detailed outline of the extension process in Step
3, on which the proof of Theorem \ref{mainthm} (2) is based.

Let $\bar f_s :F_s\to F_s$ be any extension of $\bar f_t :
e_0(V)\to e_1(V)$ (Lemma \ref{ResPA}). Let $W_1= F_s- \text{int}
e_1(V)$ and $h: W_1\to W_1$ be any pseudo-Anosov map. Let
$\mathcal A_1$ be any maximal independent set of disjoint circles
on $W_1$ (see Definition \ref{indepDef}), let $\tau (\mathcal
A_1)$ be a composition of Dehn twists along all components in
$\mathcal A_1$. Then $f_s = \tau^l (\mathcal A_1)\circ h^k \circ
\bar f_s$ is pseudo Anosov for large integers $k$ and $l$ (Lemmas
\ref{essential}, \ref{noncyclic}, \ref{indepLem}, \ref{atoroidal},
\ref{DehnTwist}).

In Section 3, we choose $\bar f_s$, $h$ and $\mathcal A_1$
carefullly so that Theorem \ref{mainthm} (2) is proved (Lemmas
\ref{HomoIden}, \ref{matrix}, \ref{5rank}).

We end the introduction by a comment on a related work \cite{WWZ}.
The main result in \cite{WWZ} is that for an orientable closed
surface $F$ with $\chi(F)<0$ and two non-isotopic circles $c$ and
$c'$ on $F$, if $g(c)={c'}$ for some automorphism $g$ on $F$, then
$f(c)= c'$ for some pseudo-Anosov map $f$ on $F$. Some arguments
in proving Lemmas \ref{noncyclic} and \ref{atoroidal} were
influenced by that in \cite{WWZ}. Indeed \cite{WWZ} is produced in
a rather earlier stage of understanding the present project.

\begin{acknowledgments}
We are grateful to Dr. Hao Zheng for drawing the figures in this
paper. M.Boileau and S.C.Wang wish to thank Prof. D.Gabai and
Prof. J.Mess for helpful conversations.

Y.Ni is partially supported by a Centennial fellowship of the
Graduate School at Princeton University. S.C.Wang is partially
supported by MSTC.

Y.Ni joined this project when he was a graduate student at Peking
University. The paper was finished when Y.Ni visited Peking
University.
\end{acknowledgments}

\section{Homotopic pinches and pseudo-Anosov extensions}

\begin{definition}\label{pinchdef} Let $D$ be a fixed disc in $F_t$ and $V=F_t-\textrm{int}(D)$.
A degree one map $p: F_s\to F_t$ is a {\it pinch} if $p|: p^{-1}
(V)\to V$ is a homeomorphism.
\end{definition}

It has been known since Nielsen and Kneser that every degree one
map between surfaces is homotopic to a pinch, see \cite{Ed} for a
reference.

\begin{notation}
In the rest of this paper, $r=s,t$ and $j=0,1$.

Recall that $M(F_r,f_r)=F_r\times[0,1]/f'_r$, where $f'_r :
F_r\times\{0\}\to F_r\times\{1\}$ is given by
$f'_r(x,0)=(f_r(x),1)$. Let $q_r: F_r\times [0,1]\to F_r$ be the
projection defined by $q_r(x,u)=x$, and $e_{r,j} : F_r\to
F_r\times \{j\}\subset F_r\times[0,1]$ be the homeomorphism given
by $e_{r,j}(x)=(x,j)$. Let $o_r:F_r\times[0,1]\to M(F_r, f_r)$ be
the quotient map and $F'_r=o_r(F_r\times 0)=o_r(F_r\times 1)$.

Then $$q_r\circ
 f'_r\circ e_{r,0} =f_r, \qquad e_{r, j}\circ q_r= \text{id}_{F_r\times \{j\}}. \eqno(*)$$
\end{notation}

\begin{lemma}\label{1fiberprsv}
There exists a fiber preserving degree one map $P: M(F_s, f_s)\to
M(F_t, f_t)$ if and only if there are homotopic pinches $p_0, p_1
: F_s\to F_t$ such that $p_1\circ f_s = f_t\circ p_0 $.
\end{lemma}

\begin{proof*}
Suppose first that $P: M(F_s, f_s) \to M(F_t, f_t)$ is a fiber
preserving degree 1 map. Up to homotopy we may assume that
$P^{-1}(F'_t)=F'_s$ and $P| : F'_s\to F'_t$ is a pinch. Moreover
we may assume that the induced degree one map on $S^1$ is
orientation preserving. Then by cutting $M(F_r, f_r)$ along
$F_r'$, $P$ provides a proper degree one map
$$\bar P:
(F_s\times[0,1], F_s\times\{0\}, F_s\times\{1\})\to
(F_t\times[0,1], F_t\times\{0\}, F_t\times\{1\})$$
 with the
property $\bar P|_{F_s\times\{1\}}\circ f'_s= f'_t \circ \bar
P|_{F_s\times\{0\}}$.

Let $p_j=q_t \circ \bar P|_{F_s\times\{j\}}\circ e_{s,j}: F_s\to
F_t$. Then $p_j$ is a pinch and $q_t\circ \bar P:
F_s\times[0,1]\to F_t$ is a homotopy from $p_0$ to $p_1$. Moreover
$\bar P|_{F_s\times\{1\}}\circ f'_s= f'_t \circ \bar
P|_{F_s\times\{0\}}$ and $(*)$ imply that
\begin{eqnarray*}
p_1\circ f_s &= &q_t \circ P|_{F_s\times\{1\}}\circ e_{s, 1}\circ
q_s\circ f'_s\circ e_{s,0} \\
&= &q_t \circ \bar P|_{F_s\times\{1\}}\circ f'_s\circ e_{s,0}\\
&= &q_t \circ f'_t \circ \bar P|_{F_s\times\{0\}}\circ e_{s,0}\\
&= &q_t \circ f'_t \circ e_{t,0}\circ q_t\circ \bar
P|_{F_s\times\{0\}}\circ e_{s,0}\\
&= &f_t\circ p_0.
\end{eqnarray*}

Suppose then there are two homotopic pinches $p_0, p_1 : F_s\to
F_t$ such that $p_1\circ f_s  = f_t\circ p_0 $. Let $P' :
F_s\times [0,1] \to F_t$ be a homotopy from $p_0$ to $p_1$.  Then
$P'$ provides a proper degree one map
$$\bar P:
(F_s\times[0,1], F_s\times\{0\}, F_s\times\{1\})\to
(F_t\times[0,1], F_t\times\{0\}, F_t\times\{1\})$$ defined by
$\bar P(x,u)=(P'(x,u),u)$. Clearly $\bar P$ is fiber preserving
and $p_j=q_t \circ \bar P|_{F_s\times\{j\}}\circ e_{s, j}$. Then
$p_1\circ f_s = f_t\circ p_0 $ and $(*)$ implies that
$$q_t \circ P|_{F_s\times\{1\}}\circ e_{s, 1}\circ
q_s\circ f'_s\circ e_{s,0} =q_t \circ f'_t \circ e_{t,0}\circ
q_t\circ \bar P|_{F_s\times\{0\}}\circ e_{s,0},
$$
hence
$$q_t \circ P|_{F_s\times\{1\}}\circ f'_s\circ e_{s,0}
=q_t \circ f'_t\circ \bar P|_{F_s\times\{0\}}\circ e_{s,0}.
$$
Since $q_t|_{F_t\times \{1\}}$ and $e_{s,0}$ are invertible, we
have

$$P|_{F_s\times\{1\}}\circ f'_s= f'_t\circ \bar P|_{F_s\times\{0\}}.$$

Hence $\bar P$ is able to induce a fiber preserving degree one map
$P: M(F_s, f_s)\to M(F_t, f_t)$.\end{proof*}

By Lemma \ref{1fiberprsv}, to prove Theorem \ref{mainthm} (1), we
need only to find two homotopic pinches $p_0, p_1: F_s\to F_t$ and
a pseudo-Anosov map $f_s:F_s\to F_s$ such that $p_1\circ f_s =
f_t\circ p_0$.

For $D\subset F_t$ and $V=F_t-\textrm{int}(D)$ given in Definition
\ref{pinchdef}, we can assume that $f_t|D=\textrm{id}$ up to
isotopy.

\begin{lemma}\label{ResPA}
If $f_t:F_t\to F_t$ is a pseudo-Anosov map, then $f_t|V:V\to V$ is
also a pseudo-Anosov map.
\end{lemma}

\begin{proof*}
Suppose there is a  non-contractible circle $c$ on $V$ such that
$f_t|^n(c)\sim c$ on $V$ for some $n>0$, then $f_t^n(c)\sim c$ on
$F_t$. Since $f_t$ is pseudo-Anosov, $c$ is contractible on $F_t$.
Hence $c$ bounds a disc $D^*$ in $F_t$ and $D\subset D^*$. It
follows that $c=\partial D^*$ is parallel to $\partial D=\partial
V$. Hence $f_t|:V\to V$ is pseudo-Anosov by definition.
\end{proof*}

Let $p_0, p_1 : F_s\to F_t$ be two pinches.  Then the pull-back of
$V$ into $F_s$ provides embeddings $e_j:V\hookrightarrow F_s$. Let
$V_j=e_j(V)$, $W_j=F_s-\textrm{int}(V_j)$, ($j=0,1$). The
following lemma is clear

\begin{lemma}\label{2extension}
$$\bar f_t := e_1 \circ f_t \circ e_0^{-1}: V_0\to V_1$$ is a
homeomorphism. Moreover, $\bar f_t$ can be extended to a
homeomorphism $f_s:F_s\to F_s$, such that $f_t \circ p_0= p_1
\circ f_s$.
\end{lemma}

A necessary condition to guarantee the extension $f_s$ in Lemma
\ref{2extension} to be pseudo-Anosov is that $e_0(\partial D)$ is
{\bf not} homotopic to $e_1(\partial D)$.

Now with Lemma \ref{1fiberprsv} and Lemma \ref{2extension},
Theorem \ref{mainthm} (1) follows from the following Lemma
\ref{3exsitpinch} and Proposition \ref{4PAextension}.

\begin{lemma}\label{3exsitpinch}
With the notation above, there exist two pinches $p_0, p_1 :
F_s\to F_t$ such that

(i) $p_0$ and $p_1$ are homotopic;

(ii) $e_0(\partial D)$ is not homotopic to $e_1(\partial D)$.
\end{lemma}

\begin{proof*}We will find two essential circles
$\gamma_0,\gamma_1\subset F_s$ such that

(1) $\gamma_0$ is not homotopic to $\gamma_1$,

(2) $\gamma_j$ separates $F_s$ into 1-punctured surfaces $V_j$ and
$W_j$, where $V_0, V_1$ have genus $g_t$.

Then we define the pinch $p_j: F_s\to F_t$  such that $W_j$ is the
pinched part.

\noindent{\bf Case 1.}  $g_t \ge 2$. $W_j$, $V_j$, are shown in
Figure~1.

Let $p_j: F_s\to F_t$ be a pinch which sends $W_j$ to $D_j\subset
F_t$ such that the restrictions $p_0|, p_1|: F_s\setminus (W_0\cup
W_1)\to F_t\setminus (D_0\cup D_1)$ are identical homeomorphisms.
This requirement can be reached if we consider $p_j: F_s\to F_t$
as a quotient map which is the identity on $F_s\setminus W_j$ and
pinches $W_j$ to $D_j$.

Note $W_0\cup W_1$ is a compact surface with two boundary
components and  $D_0\cup D_1$ is annulus. Moreover the
restrictions $p_0|, p_1|: W_0\cup W_1\to D_0\cup D_1$ are degree
one maps which are identity on the boundary, it follows from
classical argument that  $p_0|, p_1| : W_0\cup W_1\to D_0\cup D_1$
are homotopic relative to the boundary, and finally $p_0, p_1 :
F_s\to F_t$ are homotopic.

\vskip 0.5 true cm
\begin{center}
    \psfrag{0}[][]{$\gamma_0$}
    \psfrag{1}[][]{$\gamma_1$}
    \psfrag{a}[][]{$V_0$}
    \psfrag{b}[][]{$V_1$}
    \includegraphics{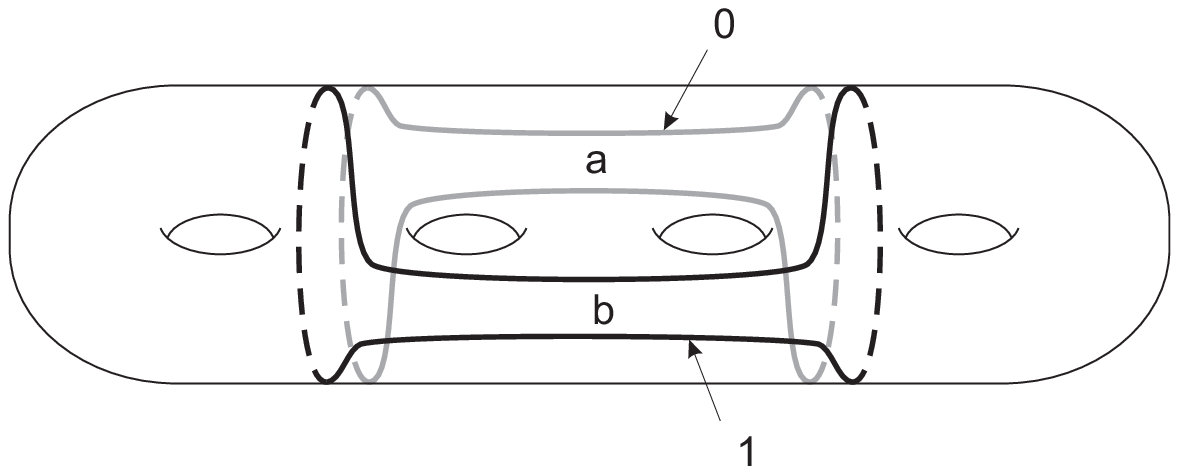}
\end{center}
\centerline {\bf Figure 1}

\vskip 0.5 true cm

\noindent{\bf Case 2.} $g_t=1$.  Then $\pi_1(F_t)=\homo_1(F_t)$ is
abelian and for each map $p:F_s\to F_t$,  $p_{\pi}: \pi_1(F_s)\to
\pi_1(F_t)$ is a composition of $\sigma: \pi_1(F_s)\to
\homo_1(F_s)$ and $p_{\#}: \homo_1(F_s)\to \homo_1(F_t)$, where
$\sigma$ is the abelianizing map, $p_{\#}$ is the map on homology.
So the homotopy class of $p$ is determined by $p_{\#}$ by
elementary homotopy theory (see \cite{H}).

Using this fact, we can construct $\gamma_0$ and $\gamma_1$ as
following: choose essential curves $\alpha$, $\beta_0$ and
$\beta_1$ on $F_s$, see Figure~2,  such that

(1) $\beta_0$ and $\beta_1$ are in the same homology class, but
not in the same homotopy class;

(2) $|\alpha\cap\beta_0|=|\alpha\cap\beta_1|=1$.

Let $\gamma_j=\partial N(\alpha\cup\beta_j)$,
$V_j=N(\alpha\cup\beta_j)$. It is easy to check that
 $\gamma_0\nsim\gamma_1$ and $p_{0\#}=p_{1\#}: \homo_1(F_s)\to
 \homo_1(F_t)$. \end{proof*}

\vskip 0.5 true cm
\begin{center}
    \psfrag{0}[][]{$\beta_0$}
    \psfrag{1}[][]{$\beta_1$}
    \psfrag{a}[][]{$\alpha$}
    \includegraphics{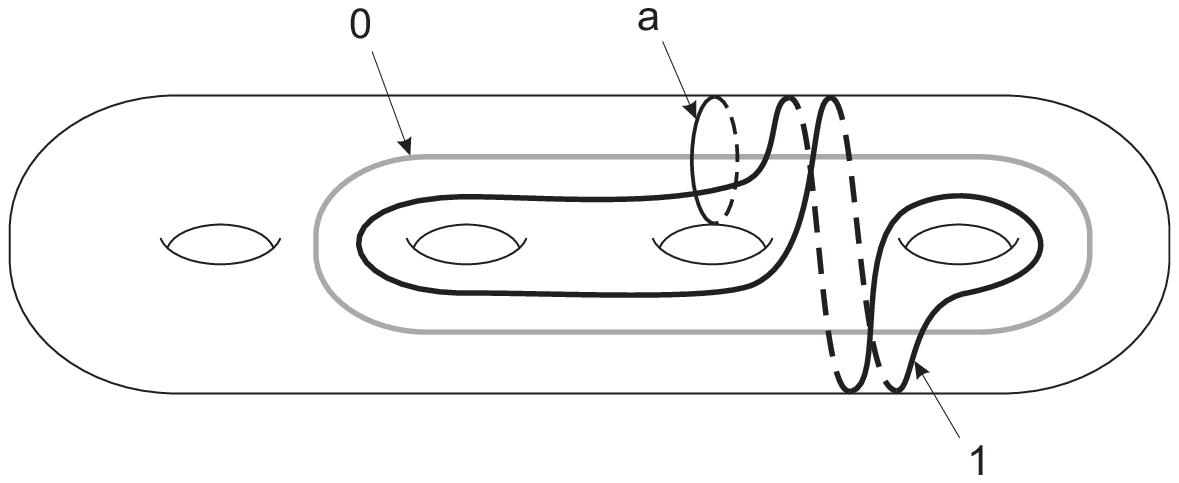}
\end{center}
\centerline {\bf Figure 2}

\vskip 0.5 true cm

\begin{proposition}\label{4PAextension}
With the notation as above, once the two pinches $p_0, p_1 :
F_s\to F_t$ are chosen to meet (i) and (ii) in Lemma
\ref{3exsitpinch}, the extension $f_s$ in Lemma \ref{2extension}
can be chosen to be pseudo-Anosov.
\end{proposition}

Suppose the two pinches $p_0, p_1 : F_s\to F_t$ are chosen to meet
(i) and (ii) in Lemma \ref{3exsitpinch}, and $f_s :(F_s,V_0)\to
(F_s,V_1)$ is an extension of $\bar f_t = e_1 \circ f_t| \circ
e_0^{-1}: V_0\to V_1$ with the condition $p_1\circ f_s = f_t\circ
p_0 $.

\begin{lemma}\label{essential}
(1) If $c$ is an essential circle in $V$, then $c$ is essential in
$F_t$.

(2) No non-trivial circle  $c\subset V_j$ can be isotoped into
$W_{j'}$, $j, j'\in\{0,1\}$, $j\ne j'$.
\end{lemma}

\begin{proof*}  (1) Otherwise $c$ would bound a disk $D^*$ in $F_t$
with $\partial V\subset D^*$, hence $c$ is parallel to $\partial
V$ in $V$, a contradiction.

(2) Otherwise say $c\subset V_0$ is a non-trivial circle, which is
isotopic to a circle $c'\subset W_1$ in $F_s$.

First suppose that $c$ is essential in $V_0$. By $p_0\sim p_1$ we
have $p_0(c)\sim p_1(c')$. On one hand $c$ is essential in $V_0$
implies that $p_0(c)$ is essential in $V_0$, and then $p_0(c)$ is
essential in $F_t$ by (1). But on the other hand, $c'\subset W_1$
implies that $p_1(c')$ is homotopically trivial. We reach a
contradiction.

Then suppose that  $\partial W_1$ can be isotoped into $W_0$. Then
one of the two components $V_1$ and $W_1$ must be contained in
$W_0$. Since $W_0$ and $W_1$ are homeomorphic, if $W_1\subset
W_0$, we must have $\partial W_0$ is parallel to $\partial W_1$, a
contradiction. Hence $V_1\subset W_0$, which implies that
$\pi_1(V_1)\subset \ker p_{0\pi}=\ker p_{1\pi}$, which clearly is
impossible.
\end{proof*}

So what remains to us is to modify $f_s|_{W_0}$. \vskip 0.3 true
cm

\begin{definition}\label{indepDef} \cite{WWZ} A set of mutually disjoint circles
$\mathcal C=\{c_1,\dots,c_m\}$
 on a compact surface $F$ is an {\it independent set}, if the circles in $\mathcal C$ are
essential and mutually non-parallel.
\end{definition}

\begin{lemma}\label{noncyclic}
Let $h:W_1\to W_1$ be a pseudo-Anosov map which is the identity in
$\partial W_1$. We extend $h$ by identity to an automorphism $h$
of $F_s$. $\mathcal A$ is a maximal independent set of circles in
$W_0$. $f=f_s:F_s\to F_s$ is an extension of $\bar f_t: V_0\to
V_1$.

Then when $k$ is sufficiently large, for any $\alpha\in \mathcal
A$, $h^kf(\alpha)$ is not isotopic to any circle in $\mathcal A$.
\end{lemma}
\begin{proof*}
Suppose $k_1<k_2$, $\alpha\in \mathcal A$. We claim that
$h^{k_1}f(\alpha)\nsim h^{k_2}f(\alpha)$ in $F_s$. In fact,
$h^{k_1}f(\alpha)$ is an essential curve in $W_1$, and any two
curves in $W_1$,which are homotopic in $F_s$, must be homotopic in
$W_1$. But $h|_{W_1}$ is a pseudo-Anosov automorphism on $W_1$, so
$h^{k_2-k_1}(h^{k_1}f(\alpha))$ is not isotopic to
$h^{k_1}f(\alpha)$.

Hence for any $\alpha\in \mathcal A$, there are only finitely many
$k$, such that $h^kf(\alpha)$ is homotopic to a circle in
$\mathcal A$. Hence the conclusion holds.
\end{proof*}

From now on we replace $f$ by $h^kf$.

Let $\mathcal A_0$ be a maximal independent set of circles on
$W_0$, $\mathcal A_1$ be its image under $f_s$. Let
$V_j'=e_{s,j}(V_j)$, $W_j'=e_{s,j}(W_j)$, $\mathcal
A_j'=e_{s,j}(\mathcal A_j)$. Let $L=o_s(\mathcal A'\bigcup\partial
W_0')$, $V'=o_s(V'_j)$, $W'=o_s(W'_j)$, $F'=o_s(F_s\times \{j\}
)$; and $X=M(F_s, f_s)-\mathrm{int}(N(L))$, $F^*=F'\cap X$. Then
$$f_s': (F_s\times
\{0\}, V_0', W_0', \mathcal A_0)\to (F_s\times \{1\},
 V_1', W_1', \mathcal A_1)$$
 is a homeomorphism. We will prove
that $X$ is hyperbolic. We first have

\begin{lemma}\label{indepLem} (1) $\mathcal A_j'$ is a maximal
independent set of circles on $W_j'$ and each non-trivial circle
in $F^*$ is either essential in $V_j'$ or parallel to a component
of $\partial W_j'\cup \mathcal A_j'$.

(2) No component of  $\partial W_0\cup \mathcal A_0$ is homotopic
to a component of $\partial W_1\cup \mathcal A_1$, in $F_s$.
\end{lemma}
\begin{proof*}  (1) follows directly from the definitions and the
constructions.

(2) $\partial W_1$ (resp. $\partial W_0$) can not be isotoped into
$W_0$ (resp. $W_1$) by Lemma \ref{essential} (2). No component of
$\mathcal A_1=f_s(\mathcal A_0)$ is isotopic to a component of
$\mathcal A_0$ by Lemma \ref{noncyclic}. Hence (2) follows.
\end{proof*}

\vskip 0.5 true cm
\begin{center}
    \psfrag{0}[][]{$V_0$}
    \psfrag{1}[][]{$V_1$}
    \psfrag{2}[][]{$\partial W_0$}
    \psfrag{3}[][]{$\partial W_1$}
    \psfrag{a}[][]{$W_0$}
    \psfrag{b}[][]{$W_1$}
    \psfrag{c}[][]{$A_0$}
    \psfrag{d}[][]{$A_1$}
    \includegraphics{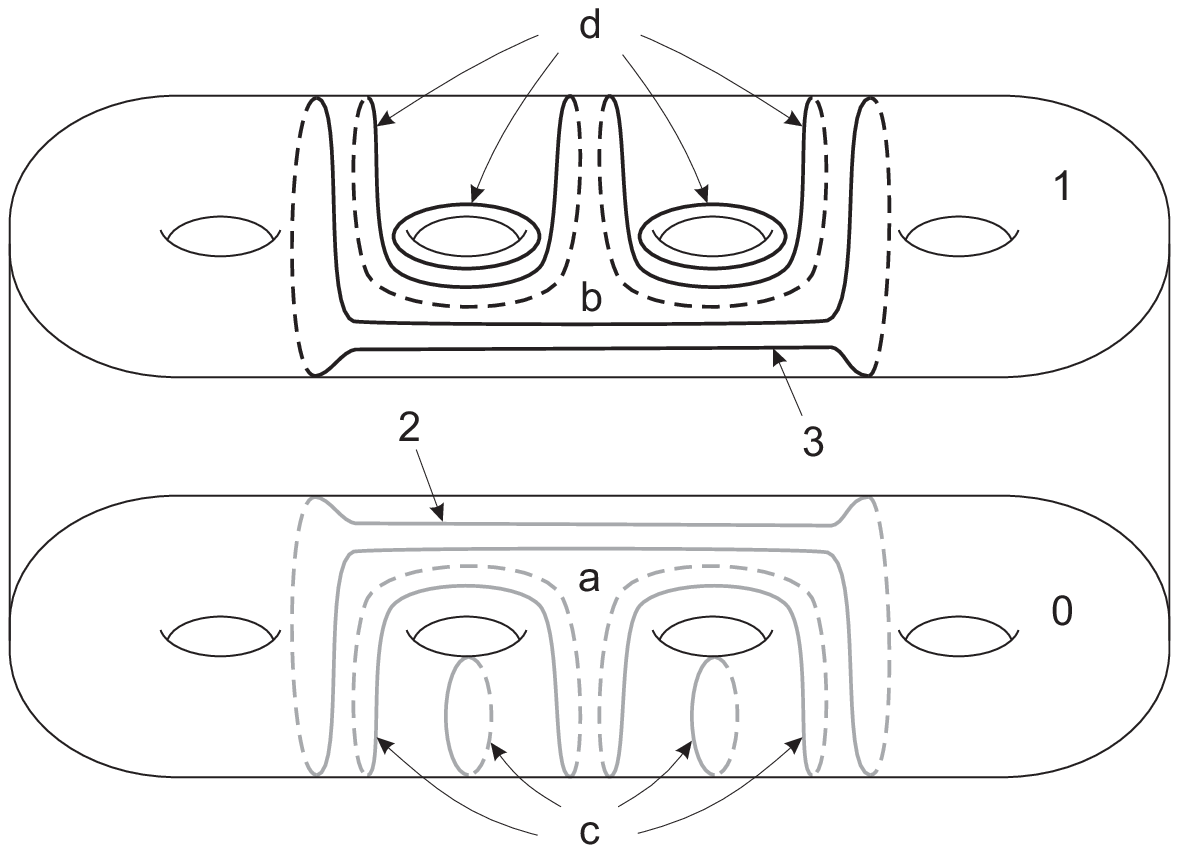}
\end{center}
\centerline{\bf Figure 3} \vskip 0.5 true cm

Figure 3 is what happens in $F_s\times [0,1]$, and clearly
illustrate the meaning of Lemma \ref{indepLem}.

\begin{lemma}\label{atoroidal}
X is atoroidal.
\end{lemma}
\begin{proof*}
Suppose that $T$ is an essential torus in $X$. We assume that $T$
has been isotoped in $X$ so that  $|T\cap F^*|$ is minimal. Then
$T\cap F^*=\mathcal C^*$ consists of $\pi_1$-injective circles on
both $T$ and $F^*$. Note $\mathcal C^*\ne\emptyset$, otherwise $T$
would be an incompressible torus in $F_s\times [0,1]$, which is
impossible.

Cutting $X$ along $F^*$, we get a manifold $X'\subset
F_s\times[0,1]$ and $T\backslash\mathcal C$ is a collection of
annuli $A_1,A_2,\dots,A_n\subset X'$. In the case $A_k$ is
vertical, we denote the component of $\partial A_k$ in
$F_s\times\{j\}$ by $c_{k,j}$, and
$q_s(c_{k,0})=q_s(c_{k,1})\subset F_s$ by $c_k$.

Now we claim that if one component of $\partial A_k$ is essential
in $V_j'$, then $A_k$ is vertical. Because otherwise $\partial
A_k\subset V_j'$, and we can push $A_k$ across $V'$ in $X$ to
reduce $|T\cap F^*|$, a contradiction.

By Lemma \ref{indepLem} (1), there are two cases:

{\bf Case 1.} Some component of $\partial A_k$ is an essential
circle on $V_j'$.

By the claim above, $A_k$ is vertical. Then $f_s'(c_{k,0})$ is a
component of $\partial A_{l}$. ($l=k-1$ or $k+1$.) Since
$f_s'(c_{k,0})$ is essential in $V_1'$, it follows by the claim
$A_l$ is vertical and $f_s'(c_{k,0})=c_{l,1}$. Then $f_s(c_k)=c_l$
and $c_l$ is essential in $V_1$. Hence $c_l$ can not be isotoped
into $W_0$ by Lemma \ref{essential} (2). Since $c_{l,0}$ is
disjoint from $\partial W'_0$, $c_l\subset V_0$. Clearly $c_l$ is
still an essential circle in $V_0$, and therefore $c_{l, 0}$ is
essential in $V'_0$. Since $T$ is connected, repeat the same
argument finitely many times, we get that

(1) all $A_k$ are vertical;

(2) all $c_{k,0}$ are essential in $V_0'$ and all $c_{k,1}$ are
essential in $V_1'$;

(3) $f_s'(c_{k, 0})=c_{k+1,1}$ (re-indexing $A_k$ if needed, and
the subscript $k$ is considered mod $n$). Hence each $c_k$ is
essential in both $V_0$ and $V_1$, and $f_s(c_k)=c_{k+1}$.

Now both $p_0(c_k)$ and $p_1(c_k)$ are essential circles in $V$,
and therefore essential in $F_t$ by Lemma \ref{essential} (1).
Since $f_t\circ p_0 =p_1\circ f_s$, we have that $f_t\circ
p_0(c_k) = p_1\circ f_s(c_k)=p_1(c_{k+1}).$ Since $p_0$ and $p_1$
are homotopic, we have $p_0(c_k)\sim p_1(c_k)$. Then up to isotopy
$f_t^{n}\circ p_0(c_k)=p_0(c_k)$, which contradicts to the fact
that $f_t$ is a pseudo-Anosov map on $F_t$.

{\bf Case 2.} Each component of $\partial A_k$ is parallel to a
component of $\partial W_0'\cup \mathcal A_0'$. By Lemma
\ref{indepLem} (2), no $A_k$ is vertical, hence both components of
$\partial A_k$ are parallel to a component $c$ of $\partial
W_0'\cup \mathcal A_0'$. So $A_k$ can be rel $\partial A_k$
isotoped into $N(c)$ in $X'$. Hence back to $X$ the torus $T$ can
be isotoped into $N(o_s(c))$. This means that $T$ is boundary
parallel in $X$, contrary to our assumption.
\end{proof*}

It is easy to see that $X$ is irreducible: a reducing sphere $S$
would bound a ball $B$ in $M(F_s,f_s)$, because $M(F_s,f_s)$, as a
surface bundle over circle, is irreducible. Hence $B$ would
contain some component of $L$. This is impossible because each
component of $L$ is essential in $M(F_s,f_s)$.

$X$ is not a Seifert fibered space: it contains $q(F_s\times
\frac12)$, a non-separating, hyperbolic, closed incompressible
surface. No such surface exists in a Seifert fibered space with
boundary, because an essential surface in such a manifold is
either horizontal (hence bounded) or vertical (hence a torus or an
annulus).

Now the geometrization theorem of Thurston for Haken manifolds
\cite{T2} leads us to the following:

\begin{corollary}\label{hyperbolic}
$X$ is a hyperbolic manifold.
\end{corollary}

Suppose $L=\{\alpha_1,\dots,\alpha_m\}$, let $T_l$ be the torus
$\partial N(\alpha_l)$ on $\partial X$, $l=1,\dots,m$. Denote by
$\tau_c$ the right hand Dehn twist along a circle $c$ on $F_s$.
Pick a meridian-longitude pair for each $T_l$, with longitude a
component of $F'\cap T_l$. $q_l$ is a slope on $T_l$, define
$X(q_1,\dots,q_m)$ to be the manifold obtained by $q_l$ Dehn
filling on $T_l$. The following lemma points a well-known relation
between Dehn fillings and Dehn twists, which has been used in some
papers, say \cite{LM} and \cite{WWZ}.

\begin{lemma}\label{DehnTwist}
Let $\tilde f
(k_1,\dots,k_m)=\tau_{\alpha_1}^{k_1}\circ\dots\circ\tau_{\alpha_m}^{k_m}\circ
f$. Then
$$M(F_s,\tilde
f(k_1,\dots,k_n))=X(1/{k_1},\dots,1/{k_m})$$ for all $k_l \in
\mathbb Z$.
\end{lemma}

\begin{proof*}[Proof of Proposition \ref{4PAextension}.] By Corollary
\ref{hyperbolic}, $X$ is a hyperbolic manifold, therefore, by the
hyperbolic surgery theorem of Thurston \cite{T1},
$X(1/{k_1},\dots,1/{k_m})$ is hyperbolic for sufficiently large
$k_l$. The previous lemma implies that
$$X(1/{k_1},\\\dots,1/{k_m})=M(F_s,\tilde
f(k_1,\dots,k_m)).$$ The theorem now follows from Thurston's
theorem that $M(F_s,\tilde f)$ is hyperbolic if and only if
$\tilde f$ is isotopic to a pseudo-Anosov map
(\cite{T2},\cite{O}).\end{proof*}

\section{Adjusting Betti numbers}
Now we pay attention to the Betti numbers of the surface bundles.
Using HHN extension one can calculate directly that
$$\homo_1(M(F,f); \mathbb Z) = \homo_1(F, \mathbb Z)/\ker (\tensor I_{2g} -
f_{\#})\:\oplus \mathbb Z$$ where $g=g(F)$, $\tensor I_{2g}$ is
the unit matrix in $SL_{2g} (\mathbb Z)$.

By abuse of notation, denote the image of $\homo_1(V_j)$ in
$\homo_1(F_s)$ by $\homo_1(\widehat{V}_j)$. Similarly, define
$\homo_1(\widehat{W}_j)$. We have
$$\homo_1(F_s)=\homo_1(\widehat{V}_0)\oplus\homo_1(\widehat{W}_0)=\homo_1(\widehat{V}_1)\oplus\homo_1(\widehat{W}_1).$$

\begin{lemma}\label{HomoIden}
$$\homo_1(\widehat{V}_0)=\homo_1(\widehat{V}_1),$$
$$\homo_1(\widehat{W}_0)=\homo_1(\widehat{W}_1),$$
as subgroups of $\homo_1(F_s)$.
\end{lemma}
\begin{proof*}
The second equation is easy, because $\homo_1(\widehat{W}_j)=\ker
p_{j\#}$, and $p_0\sim p_1$ implies $p_{0\#}=p_{1\#}$.

Now we will prove the first equation. Suppose $c_0\subset V_0$ is
a circle, then $c=p_0(c_0)$ is a circle in $V$, thus
$c_1=p_1^{-1}(c)$ is a circle in $V_1$. Since $p_0\sim p_1$, we
have a map $H:F_s\times [0,1]\to F_t,\; H(\cdot,j)=p_j(\cdot)$.

Make $H$ transverse to $c$, then $H^{-1}(c)$ is a submanifold of
$F_s\times[0,1]$ with boundary $c_0\cup c_1$. So $c_0$ and $c_1$
represent the same homology class in
$\homo_1(F_s\times[0,1])=\homo_1(F_s)$. Hence the generators of
$\homo_1(\widehat{V}_0)$ is the same as the ones of
$\homo_1(\widehat{V}_1)$, our conclusion holds.
\end{proof*}

Choose a basis of $\homo_1(\widehat{V}_0)$ and a basis of
$\homo_1(\widehat{W}_0)$ to make up a basis of $\homo_1(F_s)$.
Under this basis, $f_{s\#}$ will be represented by a matrix of the
form:
$$\begin{pmatrix}
f_{t\#} & 0\\
0 & \tensor{A}
\end{pmatrix}.$$

\begin{lemma}\label{matrix}
The map $f_s$ can be chosen so that the matrix
$\tensor{I}-\tensor{A}$ is non-degenerate.
\end{lemma}
\begin{proof*}
Let $\delta=g_s-g_t$. Choose curves
$\alpha_1,\dots,\alpha_{\delta},\beta_1,\dots,\beta_{\delta}\subset
W_1$, such that they are mutually disjoint, except that $\alpha_l$
intersects $\beta_l$ in a single point transversely. These
$2\delta$ curves form a symplectic basis of
$\homo_1(\widehat{W}_1)$. Under this basis, the intersection form
of $\homo_1(\widehat{W}_1)$ is $\begin{pmatrix}0 &
\tensor{I}_{\delta}\\-\tensor{I}_{\delta} & 0\end{pmatrix}$. So if
$f:(F_s,W_0)\to(F_s,W_1)$ is a homeomorphism,
$f_{\#}|_{\homo_1(\widehat{W}_0)}$ will be represented by a
symplectic matrix $\tensor{F}$.

We choose a map $\eta: W_1\to W_1$, such that it fixes the points
on $\partial W_1$, and it induces $\tensor{F}^{-1}$ on homology.

When $\delta>1$, by Theorem 2 in \cite{Pa}, every symplectic
matrix of rank $2\delta$ can be represented by a pseudo-Anosov map
on a closed surface of genus $\delta$. So there is a pseudo-Anosov
map $h:W_1\to W_1$, such that it fixes the points on $\partial
W_1$, and induces $-\tensor{I}_{2\delta}$ on homology. Extend
$\eta,h$ to maps on $F_s$, with the points in $V_1$ fixed.

Let $\gamma_l=\partial N(\alpha_l\cup\beta_l)$, which is a
separating circle in $W_1$ (see those separating circles in
$\mathcal A_1$, Figure 3). Now we extend
$\{\alpha_1,\dots,\alpha_{\delta},\gamma_1,\dots,\gamma_{\delta}\}$
to a maximal independent set $\mathcal{A}$ on $W_0$. Then every
curve in $\mathcal{A}-\{\alpha_1,\dots,\alpha_{\delta}\}$ is
homologous to $0$. Let $L=\mathcal A\cup\{\partial W_1\}$ So the
only Dehn twists along circles in $L$, which act nontrivially on
homology group, are
$\tau_{\alpha_1},\dots,\tau_{\alpha_{\delta}}$. The action of
products of these twists on $\homo_1(\widehat W_1)$ is represented
by a upper-triangular matrix $\tensor{T}$, whose diagonal elements
are all $1$. By Lemma \ref{noncyclic}, $h^{2k+1}\eta f$ does not
send any curve in $L$ into $L$ when $k$ is sufficiently large. The
matrix of $h^{2k+1}\eta f$, when restricted on
$\homo_1(\widehat{W}_0)$, is $-\tensor{I}_{2\delta}$. Now replace
$f$ by the composition of Dehn twists along $L$ and $h^{2k+1}\eta
f$, we have
$$\tensor{I}_{2\delta}-\tensor{A}=\tensor{I}_{2\delta}+\tensor{T}$$
is non-degenerate.

When $\delta=1$, we give a direct construction. Let $\alpha,\beta$
be a symplectic basis of $W_0$. Choose a map $h':W_1\to W_1$, with
matrix $\begin{pmatrix}2 & 1\\1 & 1\end{pmatrix}$. Then $h'\eta f$
meets the conclusion of Lemma \ref{noncyclic}. $\tau$ is the Dehn
twist along $\eta f(\alpha)$, then the matrix of $\tau$ is
$\begin{pmatrix}1 & 1\\0 & 1\end{pmatrix}$. $\tau_0$ is the Dehn
twist along $\partial W_1$. Now the matrix $\tensor A$ of
$\tau_0^m\tau^kh'\eta f$ is
$$\begin{pmatrix}1 & 1\\0 & 1\end{pmatrix}^k\begin{pmatrix}2 &
1\\1 & 1\end{pmatrix}=\begin{pmatrix}2+k & 1+k\\1 &
1\end{pmatrix}.$$ One can check $\tensor{I}-\tensor{A}$ is
non-degenerate when $k\ge0$.
\end{proof*}

To prove Theorem \ref{mainthm} (2), we need only to show the
following lemma.

\begin{lemma}\label{5rank}
With the notation as above, the extension $f_s$ in Proposition
\ref{4PAextension} can be chosen so that
$$\textrm{rank}\;( \homo_1(F_s; \mathbb Q)/\ker(\tensor I_{2g_s} -
f_{s\#}))=\textrm{rank}\; (\homo_1(F_t, \mathbb Q)/\ker(\tensor
I_{2g_t} - f_{t\#})).$$
\end{lemma}
\begin{proof*} The conclusion follows from the formula of computing
$\homo_1(M(F,f))$ and Lemma~\ref{matrix}. \end{proof*}

\begin{proof*}[Proof of Theorem \ref{mainthm}.] This theorem follows from
Proposition \ref{4PAextension} and Lemma~\ref{5rank}.
\end{proof*}

\end{document}